\documentclass[12pt]{article}
\usepackage{amsmath,amsfonts,amssymb,amsthm}
\newcounter{num}

\begin{document}
\begin{center}
{\Large\bf On some $q$-identities related to divisor functions}\\

\end{center}
\vskip 2mm \centerline{ Jiang Zeng }

\begin{center} \small  Institut Girard Desargues,
Universit\'e Claude Bernard (Lyon I)\\
21 Avenue Claude Bernard, 69622 Villeurbanne Cedex, France \\
{e-mail: \tt zeng@igd.univ-lyon1.fr}\\

\end{center}
\medskip
{\small {\bf Abstract}: We give generalizations and simple proofs
of some $q$-identities of Dilcher, Fu and Lascoux related to
divisor functions.}
\bigskip

Let $a_1,\ldots, a_N$ be $N$ indeterminates. It is  easy to see
that
\begin{equation}\label{eq:fondamental}
\frac{1}{(1-a_1z)(1-a_2z)\ldots
(1-a_Nz)}=\sum_{k=1}^N\frac{\prod_{j=1, j\neq
k}^{N}(1-a_j/a_k)^{-1}}{1-a_kz}.
\end{equation}
The coefficient of $z^\tau$ ($\tau\geq 0$) in the left side of
\eqref{eq:fondamental} is usually called the $\tau$-th
\emph{complete symmetric function}  $h_\tau(a_1,\ldots, a_N)$ of
$a_1,\ldots, a_N$. Clearly, we have $h_0(a_1,\ldots, a_N)=1$ and
equating the coefficients of $z^\tau$ ($\tau\geq 1$) in two sides
of \eqref{eq:fondamental} yields
\begin{equation}\label{eq:general}
h_\tau(a_1,\ldots, a_N):=\sum_{1\leq i_1\leq i_2\leq \cdots \leq
i_\tau\leq N}a_{i_1}a_{i_2}\ldots
a_{i_\tau}=\sum_{k=1}^N\prod_{j=1, j\neq k}^{N}(1-a_j/a_k)^{-1}
\,a_k^\tau.
\end{equation}
In particular, if $a_k=\frac{a-bq^{k+i-1}}{c-zq^{k+i-1}}$ ($1\leq
k\leq N$) for  a fixed integer $i$ ($1\leq i\leq n$), then
formula~\eqref{eq:general} with $N=n-i+1$ reads
\begin{align}
&h_\tau\left(\frac{a-bq^i}{c-zq^i},
\frac{a-bq^{i+1}}{c-zq^{i+1}},\ldots,
\frac{a-bq^n}{c-zq^n}\right)=\frac{c^{n-i+1}(zq^i/c)_{n-i+1}}{(q)_{n-i+1}(az-bc)^{n-i}}\nonumber\\
&\hspace{20pt}\cdot \sum_{k=i}^n(-1)^{k-i}{n-i+1\brack
n-k}q^{{k-i+1\choose
2}-k(n-i)}\frac{(1-q^{k-i+1})(a-bq^k)^{\tau+n-i}}{(c-zq^k)^{\tau+1}},\label{eq:master}
\end{align}
where $(x)_n=(1-x)(1-xq)\ldots (1-xq^{n-1})$ and ${n\brack
i}=(q^{n-i+1})_i/(q)_i$ with $(x)_0=1$.

The aim of this note is to show that \eqref{eq:master} turns out
to be a common source of several $q$-identities surfacing
recently in the literature.

First of all,  the $i=1$ case of formula~\eqref{eq:master} with
$\tau=m-n+1$ corresponds to an identity of Fu and
Lascoux~\cite[Prop. 2.1]{FL1}:
\begin{align}
&h_\tau\left(\frac{a-bq}{c-zq}, \frac{a-bq^{2}}{c-zq^{2}}, \ldots,
\frac{a-bq^{n}}{c-zq^{n}}\right)\nonumber\\
&=\frac{c^n(zq/c)_n}{(q)_n(az-bc)^{n-1}}\sum_{k=1}^n{n\brack
k}(-1)^{k-1}q^{{k+1\choose 2}-nk}
\frac{(1-q^k)(a-bq^k)^{m}}{(c-zq^k)^{\tau+1}}.\label{eq:fl}
\end{align}
Next, for $i=1,\ldots, n$ and $m\geq 1$ set
\begin{equation}
A_i(z):=\frac{q^i(zq)_{i-1}(q)_n}{(q)_i(zq)_{n}}h_{m-1}\left(\frac{q^i}{1-zq^i},
\ldots, \frac{q^n}{1-zq^n}\right).\label{eq:defAI}
\end{equation}
Then we have the following polynomial identity in $x$:
\begin{equation}\label{eq:J}
\sum_{k=1}^n{n\brack
k}\frac{(x-1)\cdots(x-q^{k-1})}{(1-zq^k)^m}q^{mk}=\sum_{k=1}^n(-1)^{k}{n\brack
k}\frac{q^{{k\choose 2}+mk}}{(1-zq^k)^{m}}+\sum_{i=1}^nA_i(z)x^i.
\end{equation}
Indeed, using the $q$-binomial formula \cite[p. 36]{An}:
$$
(x-1)(x-q)\ldots (x-q^{N-1})=\sum_{j=0}^N{N\brack
j}(-1)^{N-j}x^jq^{N-j\choose 2},
$$
we see that the coefficient of $x^i$ ($1\leq i\leq n$) in the left
side of \eqref{eq:J}  is equal to
\begin{equation}\label{eq:AI}
\sum_{k=i}^n(-1)^{k-i}{n\brack k}{k\brack
i}\frac{q^{mk+{k-i\choose 2}}}{(1-zq^k)^{m}}
=\frac{q^i(zq)_{i-1}(q)_n}{(q)_i(zq)_{n}}h_{m-1}\left(\frac{q^i}{1-zq^i},
\ldots, \frac{q^n}{1-zq^n}\right),
\end{equation}
where the last equality follows from  \eqref{eq:master} with
$a=0$, $c=1$, $b=-1$ and $\tau=m-1$.

Now, with $z=i=1$ and  $m$ shifted to $m+1$,  formula
\eqref{eq:AI} reduces to Dilcher's identity~\cite{Dil}:
\begin{equation*}\label{eq:dilcher}
\sum_{k=1}^n{n\brack k}\frac{(-1)^{k-1}q^{{k\choose
2}+mk}}{(1-q^k)^m}=h_{m}\left(\frac{q}{1-q}, \ldots,
\frac{q^n}{1-q^n}\right)=\sum_{i=1}^nA_i(1).
\end{equation*}
On the other hand,  formula \eqref{eq:fondamental} with  $N=n+1$
and $a_i=q^{i-1}$ ($1\leq i\leq N$)  yields
$$
 \sum_{k=0}^n(-1)^k{n\brack
k}\frac{q^{{k\choose 2}+k}}{1-zq^k}=\frac{(q)_n}{(z)_{n+1}}.
$$
Hence, setting, respectively, $z=1$ and $m=1$ in formula
\eqref{eq:J}
 we recover  two recent  formulae of  Fu and
Lascoux~\cite{FL2} (see also \cite{Pr}):
\begin{equation}\label{eq:flgeneral1}
\sum_{k=1}^n{n\brack
k}\frac{(x-1)\cdots(x-q^{k-1})}{(1-q^k)^m}q^{mk}=\sum_{i=1}^n(x^i-1)A_i(1),
\end{equation}
and
\begin{equation}\label{eq:flgeneral2}
\sum_{k=0}^n{n\brack k}\frac{(x-1)\cdots(x-q^{k-1})}{1-zq^k}q^{k}=
\frac{(q)_n}{(z)_{n+1}}\sum_{i=0}^n\frac{(z)_i}{(q)_i}x^iq^i.
\end{equation}

\section*{Acknowledgement} The author is
partially supported by EC's IHRP Programme, within the Research
Training Network ``Algebraic Combinatorics in Europe'', grant
HPRN-CT-2001-00272.

\end{document}